\documentclass[12pt,a4paper]{amsart}
\usepackage{amsfonts}
\usepackage{amsthm}
\usepackage{amsmath}
\usepackage{amscd}
\usepackage[utf8]{inputenc}
\usepackage{t1enc}
\usepackage[mathscr]{eucal}
\usepackage{indentfirst}
\usepackage{graphicx}
\usepackage{graphics}
\usepackage{pict2e}
\usepackage{epic}
\numberwithin{equation}{section}
\usepackage[margin=2.9cm]{geometry}
\usepackage{epstopdf}

\theoremstyle{plain}
\newtheorem{Th}{Theorem}[section]
\newtheorem{Lemma}[Th]{Lemma}

 \theoremstyle{definition}
\newtheorem{Def}[Th]{Definition}

\newtheorem{?}[Th]{Problem}

\begin{document}

\title{THE UNCERTAINTY PRINCIPLE FOR THE TWO-SIDED QUATERNION FOURIER TRANSFORM}

\author[Y. El Haoui]{Youssef El Haoui}

\address{Youssef El Haoui\\ Department of Mathematics and Computer Sciences, Faculty of Sciences, Equipe d'Analyse Harmonique et Probabilit\'es, University Moulay Ismail, BP 11201 Zitoune, Meknes, Morocco}
\email{youssefelhaoui@gmail.com}

\author[S. Fahlaoui]{Said Fahlaoui}

\address{Said Fahlaoui \\ Department of Mathematics and Computer Sciences, Faculty of Sciences, Equipe d'Analyse Harmonique et Probabilit\'es, University Moulay Ismail, BP 11201 Zitoune, Meknes, Morocco}
\email{saidfahlaoui@gmail.com}

 \subjclass[2010]{}


\begin{abstract} In this paper, we provide the Heisenberg's inequality and the Hardy's theorem for the two-sided quaternion Fourier transform.
\end{abstract}

\maketitle

\section{Introduction} In harmonic analysis, the uncertainty principle states that a non zero function and its Fourier transform cannot both be very rapidly decreasing. This fact is expressed by two formulations of the uncertainty principle for the Fourier transform, Heisenberg's inequality and Hardy's theorem.\

The classical Fourier transform, defined over the real line, is given by 

\[{\mathcal F}\left\{f\left(x\right)\right\}\left(y\right)=\int^{+\infty }_{-\infty }{e^{-2\pi ixy}}\ f\left(x\right)dx,\] 

For $f \in L^1(\mathbb{R}).$ 

The Heisenberg's inequality asserts that for $f \in L^2(\mathbb{R})$

($\int_{\mathbb{R}}{{\left|x\right|}^2}{\left|f(x)\right|}^2dx$) ($\int_{\mathbb{R}}{{\left|y\right|}^2}{\left|{\mathcal F}\left\{f\left(x\right)\right\}\left(y\right)\right|}^2dy$)$\ge \frac{{\left\|f\right\|}^4_{L^2{  (}{\mathbb{R}}{  )}}}{16{\pi }^2},$

With equality only if $f$ almost everywhere equal to  $\alpha e^{-\beta x^2}$\textit{,}$\ $ for some $\beta >0$. The proof of this inequality is given by Weyl in [11].

Hardy's theorem [7] says : If we suppose $\alpha$ and $\beta$ be positive constants and $f$ a function on the real line satisfying $\left|{  f(x)}\right|{  <}{{  Ce}}^{{  -}\alpha {{  x}}^{{  2}}}$  and $\left|{  \mathcal F}\left\{f\right\}\left({  y}\right)\right|{  <}C{{  e}}^{{  -}\beta {{  y}}^{{  2}}}$

for some positive constant $C$, then (i) $f{  =0\ }$  if $\alpha \beta >{\pi }^2$; (ii)$f{ =A}{{ e}}^{{ -}\alpha {{ x}}^{{ 2}}}$ if  $\alpha \beta {  =}{\pi }^2$; (iii) there are many $f\ $if $\alpha \beta <{\pi }^2$.

In this paper, we provide an analogues of the the Heisenberg's inequality and the Hardy's theorem for the two-sided quaternion Fourier transform.

Our paper is organized as follows. In section 2, we review  basic notions and notations related to the quaternion algebra. In section 3, we recall the definition and some results  for  the two-sided quaternion Fourier transform useful in the sequel. In section 4, we prove the Heisenberg's inequality and the Hardy's theorem for the two-sided quaternion Fourier transform. In section 5 , we provide the Hardy's theorem for the two-sided quaternion Fourier transform.

Note that we will often use the shorthand  $x{{:=}}(x_1,x_2)$,$\ y{{:=}}(y_1,y_2).$

\section{The algebra of quaternions }
The quaternion algebra over $\mathbb{R}$, denoted by $\mathbb{H}$, is an associative noncommutative four-dimensional algebra, it was invented by  W. R. Hamilton in 1843.\\
$\mathbb{H}=\{q=q_0+iq_1+jq_2+kq_3;\ q_0,q_1,\ q_2,q_3 \in \mathbb{R}$\}
where ${  i},{  \ }{  j},{  \ }{  k}$ satisfy Hamilton's multiplication rules

$ij=-ji=k~ ; jk=-kj=i~ ; ki=-ik=j~ ; i^2 = j^2 = k^2=-1.$

Quaternions are isomorphic to the Clifford algebra ${Cl}_{(0,2)}$ of ${\mathbb R}^{(0,2)}$:

$\mathbb{H} \cong {Cl}_{(0,2)}$          \hfill(2.1)

The scalar part of  a quaternion  $q \in  \mathbb{H}$\ is $q_0$ denoted by $Sc(q)$, the non scalar part(or pure quaternion) of $q$\ is $iq_1+jq_2+kq_3$ denoted by $Vec(q)$.

We define the conjugation of $q \in \mathbb{H}$\ by ~:

$\overline{q}$=$q_0-iq_1-jq_2-kq_3$. \\ The quaternion conjugation is a linear anti-involution 

\[\overline{qp}= \overline{p} \ \overline{q} ,\ \overline{p+q}= \overline{p}+\overline{q},\ \overline{\overline{p}}=p.\] 

The modulus of a quaternion  q is defined by:$\ $

\[{|q|}_Q=\sqrt{q\overline{q}}=\sqrt{{q_0}^2{  +}\ {q_1}^2{  +}{q_2}^2+{q_3}^2}.\] 

In particular, when $q=q_0$ is a real number, the module ${|q|}_Q$ reduces to the ordinary Euclidean module $\left|q\right|=\sqrt{{q_0}^2}$.

we have~:

\[{|pq|}_Q={|p|}_Q{|q|}_Q.\] 

It is easy to verify that  0$\ne q \in \mathbb{H}$ implies~:

\[q^{-1}=\frac{\overline{q}}{{{|q|}_Q}^2}.\] 

Any quaternion  $q$ can be written as $q$=\ ${{|q|}_Qe}^{\mu \theta }$  where $e^{\mu \theta }$ is understood in accordance with Euler's formula $e^{\mu \theta }={\cos  \left(\theta \right)\ }+\mu \ {\sin  \left(\theta \right)\ }$ where 
$\theta $=\ artan $\frac{\left|Vec\left(q\right)\right|_Q}{Sc\left(q\right)}$,\ 0$\le \theta \le \pi $ and \ $\mu $ :=\ $\frac{Vec\left(q\right)}{\left|{  Vec}\left({  q}\right)\right|_Q}$ verifying ${\mu }^2 =\ -1$.

In this paper, we will study the  quaternion-valued signal $f:{\mathbb{R}}^2\to \mathbb{H} $ , $f$ which can be expressed as $$f=f_0+i f_1+jf_2+kf_3,$$ with $f_m~: {\mathbb R^2}\to \ {\mathbb R}\ for\ m=0,1,2,3.$

Let the inner product of  ~$f,g\in L(\mathbb{R}^2, \mathbb{H})$ be defined by

$<f\left(x\right),g\left(x\right)>:=\int_{{\mathbb R}^2}{f}$(x)$\ \overline{g(x)}dx$

we define ${|f|}^2_{2,Q}:=~<f,f>$\ and ${|f|}^2_{1,Q}= \int_{{\mathbb R}^2}{{|f(x)|}_Q}\ dx$.

\section{The two-sided quaternion Fourier transform}
Ell[6] defined the quaternion Fourier transform (QFT) which has an important role in the representations of signals due to transforming a real 2D signal into a quaternion-valued frequency domain signal,  QFT belongs to the family of Clifford Fourier transformations because of (2.1).

There are three different types of QFT, the left-sided QFT , the right-sided QFT , and two-sided QFT [8].

We now review the definition and some properties of the two-sided QFT.
\begin{Def}
Let$\ f$ in    $L\left({\mathbb R}^2,{\mathbb H}\right)$. Then two-sided quarternionic Fourier transform of the function  $f$ is given by

${\mathcal F}\{f(x)\}$($\xi $)=$\int_{{\mathbb R}^2}{e^{-i{2\pi \xi }_1x_1}} f$($x$)$ e^{-j{2\pi \xi }_2x_2}dx,\ dx=dx_1dx_2$ 

Where ~$\xi,x\in {\mathbb R}^2.$
\end{Def}
We define a new module of $\mathcal F\{f\} $ as follows :

$${\left\|\mathcal F\left\{f\right\}\right\|}_Q\ :=\ \sqrt{\sum^{m=3}_{m=0}{{\left|\mathcal F\left\{f_m\right\}\right|}^2_Q}}. $$ 

Furthermore, we define a new $L^2$-norm of ${\mathcal F\{f\}}$ as follows:$$ {\left\|\mathcal F\{f\}\right\|}_{2,Q}:=\sqrt{\int_{\mathbb R^2}{{\left\|\mathcal F\left\{f\right\}(y)\right\|}^2_Qdy}}. $$

\begin{Lemma}{Inverse QFT} [1, Thm 2.5]

If  $f,{\mathcal F}\{f\}\in L\left({\mathbb{R}}^2,{\mathbb{H}}\right)$, then

$f(x)=\int_{{\mathbb{R}}^2}{e^{i2\pi {\xi }_1x_1}} {\mathcal F}\{f(x)\}(\xi ) e^{j{2\pi \xi }_2x_2}d\xi,\ d\xi=d\xi_1d\xi_2. $
\end{Lemma}

\begin{Lemma}{Plancherel theorem for QFT} [2, Thm. 3.2]

If $f\in $ $L^2\left({\mathbb{R}}^2,{\mathbb H}\right)~$, then 
                               
${|f|}_{2,Q}={\left\|{\mathcal F}\left\{f\right\}\right\|}_{2,Q}.$

\end{Lemma}
\begin{Lemma}{Derivative theorem (QFT)} [1, Thm 2.10] 

 If $f, \frac{{\partial }^{m+n}}{{\partial }^m_{x_1}{\partial }^n_{x_2}}f\in $ $L^2\left({\mathbb{R}}^2,{\mathbb H}\right)$  for $m,n\in {\mathbb N}$

then $\mathcal F$\{$\frac{{\partial }^{m+n}}{{\partial }^m_{x_1}{\partial }^n_{x_2}}f(x)$\}$\left(\xi \right){  =}{(2\pi )}^{m+n}{(i{\xi }_1)}^m{\mathcal F}$\{$f(x)$\}$\left(\xi \right){(j{\xi }_2)}^n.$

\end{Lemma}

\begin{Lemma}
 ${\mathcal F}\left\{e^{-\pi {\left|x\right|}^2}\right\}\left(y\right)$=$e^{-\pi {\left|y\right|}^2}$
where ~$x,y \in {\mathbb R}^2.$
 \end{Lemma}
Proof.${\mathcal F}\left\{e^{-\pi {\left|x\right|}^2}\right\}\left(y\right)=\int_{{{ \mathbb R}}^2}{e^{-2\pi ix_1y_1 }e^{-\pi (x^2_1+x^2_2)}e^{-2\pi jx_2y_2}dx_1dx_2}$ 

\hspace*{3.9cm}=$e^{{ -}\pi {(y}^2_1+y^2_2)}\int_{{ \mathbb R}}{e^{-\pi  {(x_1+{i}y_1)}^2 }dx_1}\int_{{ \mathbb R}}{e^{-\pi  {(x_2+ { j}y_2)}^2 }dx_2}.$ (Fubini)

We know that  $\int_{{ \mathbb R}}{e^{-{ z} {(t+z')}^2 }dt=}\sqrt{\frac{\pi }{z}}$,  for  ${ z},{ z'}\in {\mathbb C}$,\ Re(z)$>$0(Gaussian integral

with complex offset) 

Therefore $\int_{{ \mathbb R}}{e^{-\pi  {(x_1+i y_1)}^2 }dx_1}=\int_{{ \mathbb R}}{e^{-\pi  {(x_2+ j y_2)}^2 }dx_2}$=1 which give us the desired result.

\begin{Lemma}

The two-sided quaternion Fourier transform  ${\mathcal  F}$ maps ${{  L}}^{{  1}}\left({\mathbb R}^2,{\mathbb H}\right)$  into ${{  C}}_0\left({{ \mathbb R}}^{{  2}},{\mathbb H}\right){  \ }$and it is one-to-one, where $C_0\left({\mathbb R}^2,{\mathbb H}\right)$ is the set  of all functions ${  f}{  \ }= f_0  +i f_1+j f_2 + k f_3$  such that ${{  f}}_r$ is a real valued continuous function vanishing at infinity,  for r=0,1,2,3.

\end{Lemma}
\ \  Proof  is obtained by combining the Riemann-Lebesgue Lemma, the continuity of  ${\mathcal  F}{  \{  f}{  \}}$  for  $f$ in $L^{{  1}}\left({\mathbb R}^2,\mathbb H\right)$  [2, Thm 3.1 and Proposition 3.1(ix)] and lemma 3.2.

\begin{Lemma}
Let $f\in $ $L^1\cap L^2$($\mathbb R{^2} ,\mathbb{H})\ $. If   $\frac{\partial }{\partial x_k}$ $f$ exists, and are in $L^2$($\mathbb R{^2},\mathbb{H})$ for k=1,2, then

\[{(2\pi )}^2\int_{\mathbb R{^2}}{{\xi }^2_k}{\left\|\ {\mathcal{F}}\left\{f\left(x\right)\right\}(\xi )\ \ \right\|}^2_Q d\xi =\int_{\mathbb R{^2}}{{|\frac{\partial }{\partial x_k} f\left(x\right) |}^2_Q\ dx}.\] 
\\
Proof. For k=1\\
$\int_{\mathbb R{^2}}{{| \frac{\partial }{\partial x_1} f\left(x\right)|}^2_Q\ dx}=\int_{\mathbb R{^2}}{{|i^{-1}\ \frac{\partial }{\partial x_1} f\left(x\right)|}^2_Q\ dx}$

\hspace*{3.5 cm}=$\int_{\mathbb R{^2}}{{\left\|\ {\mathcal{F}}\left\{i^{-1}\frac{\partial }{\partial x_1} f\left(x\right)\ \right\}(\xi )  \right\|}^2_Q}$ $d\xi $  \           (lemma 3.3)

\hspace*{3.5 cm}=${(2\pi )}^2\int_{\mathbb R{^2}}{{{\xi }^2_1\left\|\ {\mathcal{F}}\left\{ f\left(x\right)\ \right\}(\xi ) \right\|}^2_Q}$ $d\xi. $                   \    (lemma 3.4)

Similar proof for  k=2.
\end{Lemma}

\begin{Lemma}

for $f\left(x\right)$=${|f\left(x\right)|}_Qe^{u\left(x\right)\theta (x)},\ $

If $\frac{\partial }{\partial x_k}$ $f$ exists for k=1,2, then

\[{|\frac{\partial }{\partial x_k}f\left(x\right)|}^2_Q={(\frac{\partial }{\partial x_k} {|f\left(x\right)|}_Q)}^2+{\ {|f\left(x\right)|}^2_Q\ \ |\frac{\partial }{\partial x_k}e^{u\left(x\right)\theta (x)}\ |}^2_Q.\] 
\end{Lemma}

Proof.

For  $f\left(x\right)$=${|f\left(x\right)|}_Qe^{u\left(x\right)\theta (x)},\ $we have~:

$\frac{\partial }{\partial x_k}f\left(x\right)=\frac{\partial }{\partial x_k}({|f\left(x\right)|}_Qe^{u\left(x\right)\theta (x)})$

\hspace*{1.5cm}=($\frac{\partial }{\partial x_k}{|f\left(x\right)|}_Q$)$\ e^{u\left(x\right)\theta (x)}$+${|f\left(x\right)|}_Q\frac{\partial }{\partial x_k}e^{u\left(x\right)\theta (x)}.$

thus

${|\frac{\partial }{\partial x_k}f\left(x\right)|}^2_Q=\frac{\partial }{\partial x_k}f\left(x\right)\overline{\frac{\partial }{\partial x_k}f\left(x\right)}$

=[($\frac{\partial }{\partial x_k}{|f\left(x\right)|}_Q$)$\ e^{u\left(x\right)\theta (x)}$+${|f\left(x\right)|}_Q\frac{\partial }{\partial x_k}e^{u\left(x\right)\theta (x)}
][( \frac{\partial }{\partial x_k}{|f\left(x\right)|}_Q$)$\ e^{-u\left(x\right)\theta (x)}$+${|f\left(x\right)|}_Q$ $\frac{\partial }{\partial x_k}e^{-u\left(x\right)\theta (x)}]$\\
=${(\frac{\partial }{\partial x_k}{|f\left(x\right)|}_Q)}^2$+${|f\left(x\right)|}^2_Q\frac{\partial }{\partial x_k}(e^{u\left(x\right)\theta (x)})\overline{\frac{\partial }{\partial x_k}\ e^{u\left(x\right)\theta (x)}}
{ +}{|f\left(x\right)|}_Q(\frac{\partial }{\partial x_k}{|f\left(x\right)|}_Q)$\\
\hspace*{1 cm}$[e^{u\left(x\right)\theta (x)}(\frac{\partial }{\partial x_k}e^{-u\left(x\right)\theta \left(x\right)}){ +}(\frac{\partial }{\partial x_k}e^{u\left(x\right)\theta \left(x\right)}) e^{-u\left(x\right)\theta \left(x\right)})].$\\
 \hspace*{12 cm}($\overline{\frac{\partial }{\partial x_k}\ e^{u\left(x\right)\theta (x)}}$=$\ \frac{\partial }{\partial x_k}\ e^{-u\left(x\right)\theta (x)}$)

where

$\left(\ \frac{\partial }{\partial x_k}e^{u\left(x\right)\theta \left(x\right)}\right)e^{-u\left(x\right)\theta \left(x\right)}+e^{u\left(x\right)\theta \left(x\right)}\left(\ \frac{\partial }{\partial x_k}e^{-u\left(x\right)\theta \left(x\right)}\right)$ =$\frac{\partial }{\partial x_k}\left(e^{u\left(x\right)\theta \left(x\right)}e^{-u\left(x\right)\theta \left(x\right)}\right)$

\hspace*{9 cm}=$\frac{\partial }{\partial x_k}1=0,$

and  $\frac{\partial }{\partial x_k}e^{u\left(x\right)\theta \left(x\right)}\overline{\frac{\partial }{\partial x_k}\ e^{u\left(x\right)\theta (x)}}={\ |\frac{\partial }{\partial x_k}e^{u\left(x\right)\theta (x)}\ |}^2_Q.$

This completes the proof.\\ \\
  
From  lemmas 3.7 and 3.8 the following theorem follows

\begin{Th}

Let $f\left(x\right)$=${|f\left(x\right)|}_Qe^{u\left(x\right)\theta (x)},\ $

if $f \in $ $L^1\cap L^2$($\mathbb R{^2} ,\mathbb{H})\ $and    $\frac{\partial }{\partial x_k}$ $f$ exists and is in $L^2$($\mathbb R{^2},\mathbb H)$, for k=1,2, then

${(2\pi )}^2\int_{\mathbb R{^2}}{{\xi }^2_k}{\left\| {\mathcal{F}}\left\{f\left(x\right)\right\}(\xi ) \right\|}^2_Q d\xi = \int_{\mathbb R{^2}}{{(\frac{\partial }{\partial x_k} {|f\left(x\right)|}_Q)}^2\ dx}+\int_{\mathbb R{^2}}
 { {|f\left(x\right)|}^2_Q \ |\left( \frac{\partial }{\partial x_k}e^{u\left(x\right)\theta (x)}\right)|}^2_Q\ dx. $ 

\end{Th}

\section{Heisenberg's inequality }
\begin{Th}

Let $f\left(x\right)={|f\left(x\right)|}_Qe^{u\left(x\right)\theta (x)}\ $.If $f,\ \frac{\partial }{\partial x_k}$ $f,\ x_k$ $f\in $ $L^2$($\mathbb R{^2} ,\mathbb{H})\ $for k=1,2   then

${|x_kf\left(x\right)|}^2_{2,Q}\ {\left\|{\xi }_k{\mathcal{F}}\left\{f\left(x\right)\right\}(\xi )\right\|}^2_{2,Q}\ \ge \frac{1}{{16\pi }^2}{|f\left(x\right)|}^4_{2,Q}$+${COV}^2_{x_k},$ 

with ${COV}_{x_k}:=\ \frac{1}{2\pi }\ \int_{\mathbb R{^2}}{{ {|f\left(x\right)|}}^2_Q~~\ ~{|x_k\left(\ \frac{\partial }{\partial x_k}e^{u\left(x\right)\theta (x)}\right)|}_Q\ dx}.$  

The equation holds if and only if 

$f\left(x\right)=De^{-a_kx^2_k}e^{u\left(x\right)\theta (x)} $ and  $\frac{\partial }{\partial x_k}e^{u\left(x\right)\theta (x)}$=$b_kx_k$\ where ${a_1,a}_2$$>$0, $D\in $  $\mathbb R{^+}$ and ${b_1,b}_2\ $are pure quaternions.

\end{Th}

Proof.

Since $\frac{\partial }{\partial x_k}f\in L^2(\mathbb R^2,\mathbb H)$, applying Lemma(3.4), we obtain that ${\xi }_1{  \mathcal F}\left\{f\left(x\right)\right\}(\xi )\in L^2$($\mathbb R^2 ,\mathbb H)$ and hence ${\left\|{  \mathcal F}\left\{f\left(x\right)\right\}(\xi )\right\|}^2_{1,Q}=\int_{\mathbb R}{\int_{\mathbb R}{\ \frac{(1+|{\xi }_1|)}{(1+|{\xi }_1|)}}{\left\|{  \mathcal F}\left\{f\left(x\right)\right\}(\xi )\right\|}_{{  Q}}{d\xi }_1{d\xi }_2}$
\\ 
\hspace*{6.5 cm}
$\le \int_{\mathbb R}{\frac{{d\xi }_1}{{(1+|{\xi }_1|)}^2}\int_{\mathbb R}{\ {(1+|{\xi }_1|)}^2}{\left\|{\mathcal  F}\left\{f\left(x\right)\right\}(\xi )\right\|}^{{  2}}_{{  Q}}{d\xi }_2}$$<$+$\ \infty $\\ \hspace*{11 cm}(by Holder's inequality)

and hence ${  \mathcal F}\left\{f\left(x\right)\right\}\in L^1$($\mathbb R^2\ ,\mathbb H)$. \\Since\ ${{ f(x)=\mathcal F}\left\{{  \mathcal F}\left\{{  f}\right\}\right\}\left(-{  x}  
  \right) \forall x \in {\mathbb R}^2}$, we have by lemma(3.6)$ \
f \ \in {  C}_0\left({\mathbb R}^2,{\mathbb H}\right)$. Thus, if \textbar $x_k\textbar \to +\infty $, then ${  f(x)}\to $0. Using the same argument as $x_k$ $f,\ \frac{\partial }{\partial x_k}$ ${(x}_kf),\ \in $ $L^2$($\mathbb R^2\ ,\mathbb H)$, we get that $x_kf\to 0$ as \textbar $x_k$\textbar $\to $+$\infty $. We get consequently $x_kf^2\to 0 \  as \ x_k\to +\infty $.                                                           
\hfill(4.1)

Now we have according to (Thm 3.9)

${(2\pi )}^2{|x_kf\left(x\right)|}^2_{2,Q})({\left\|{\xi }_k{\mathcal{F}}\left\{f\left(x\right)\right\}(\xi )\right\|}^2_{2,Q}={|x_kf\left(x\right)|}^2_{2,Q}\int_{\mathbb R{^2}}{{(\frac{\partial }{\partial x_k} {|f\left(x\right)|}_Q)}^2\ dx}$\\ \hspace*{8 cm}+${|x_kf\left(x\right)|}^2_{2,Q}\int_{\mathbb R{^2}}{ \ {|f\left(x\right)|}^2_Q \ |\frac{\partial }{\partial x_k}e^{u\left(x\right)\theta (x)}|}^2_Q\ dx. $

However

${|x_kf\left(x\right)|}^2_{2,Q}\ {|\frac{\partial }{\partial x_k} {|f\left(x\right)|}_Q|}^2_{2,Q}$$\ge $ ${(\int_{\mathbb R{^2}}{{|x_kf\left(x\right)\frac{\partial }{\partial x_k} {|f\left(x\right)|}_Q|}_Q}dx)}^2$           \hfill(4.2)

\hspace*{4.5 cm}${(=|x_kf\left(x\right)\frac{\partial }{\partial x_k} {|f\left(x\right)|}_Q|}^2_1$~~(H\^older))

\hspace*{4.5 cm}$\ge \ {(\int_{\mathbb R{^2}}{x_k{|f\left(x\right)|}_Q\frac{\partial }{\partial x_k} {|f\left(x\right)|}_Q}dx)}^2$

\hspace*{4.5 cm}=[$\frac{1}{2}{(\int_{\mathbb R{^2}}{\frac{\partial }{\partial x_k}{(}x_k{|f\left(x\right)|}^2_Q})dx-\int_{\mathbb R{^2}}{{|f\left(x\right)|}^2_Q}dx)]}^2$

\hspace*{4.5 cm}=$\frac{1}{4}{|f\left(x\right)|}^4_{2,Q}.$            \ \ (by 4.1)      

On the other hand,

${|x_kf\left(x\right)|}^2_{2,Q}\int_{\mathbb R{^2}} {\ {|f\left(x\right)|}^2_Q\ \ |\frac{\partial }{\partial x_k}e^{u\left(x\right)\theta (x)} |}^2_Q\ dx $

\hspace*{2 cm}$\ge \ \int_{\mathbb R{^2}}{{ {|f\left(x\right)|}_Q}^2{|x_k\left(\ \frac{\partial }{\partial x_k}e^{u\left(x\right)\theta (x)}\right)|}_Q\ dx}$                (H\^older)\hfill(4.3)

\hspace*{2 cm}=\ ${(2\pi )}^2{COV}^2_{x_k}.$

The equation holds in    (4.2)  if and only if
$\frac{\partial }{\partial x_k} {|f\left(x\right)|}_Q$=${{{\ \alpha }_k\ x}_k|f\left(x\right)|}_Q$  (${\alpha }_k\in $ $\mathbb{R}$)

which can be written as $\frac{\partial }{\partial x_k} {|f\left(x\right)|}_Q$=$-{{\  2}a_k \ x_k|f\left(x\right)|}_Q$  ($a_k\in  \mathbb{R}$).

The solution is the Gaussian~:

${|f\left(x\right)|}_Q=De^{-a_kx^{{\  2}}_k}$,\ $D\in \mathbb{R}^+$ ,  
$a_k >0$  to keep $f\in L^2$($\mathbb R{^2} ,\mathbb{H})$.

The equation holds in    (4.3)  if and only if

$\frac{\partial }{\partial x_k}e^{u\left(x\right)\theta (x)}$=$b_kx_k$ where ${b_1,b}_2\ $are pure quaternions.

\section{Hardy's theorem }
In this section, we prove Hardy's theorem for the two-sided quaternion Fourier transform.

We introduce a basis${{\{\varphi }_{k,l}\}}_{k,l\in {\mathbb N}}$ of ${\mathcal S}({\mathbb R}^2,{\mathbb H})$   and  $L^2({\mathbb R}^2,{\mathbb H})$   which is defined as follows :

With ${\varphi }_{k,l}\left({{x}_1,x}_2\right):={\varphi }_k\left({ x}_1\right){\varphi }_l\left({ x}_2\right)$\  for $\left(x_1,x_2\right)\in {\mathbb R}^2,$

And  ${\varphi }_k\left(x\right)=\frac{{(-1)}^k}{k!}{{  e}}^{\pi x^2}\frac{d^k}{{dx}^k}$(${{  e}}^{{  -}{  2}\pi x^2}).$

This basis is the tensor product basis, and every function $h$ in $L^2({\mathbb R}^2,{\mathbb H})$ can be decomposed as \[h=\sum^{+\infty }_{k,l=0}{a_{k,l}}\ {\varphi }_{k,l}\ \ with \ a_{k,l}\in {\mathbb H}.\]

We take note that, although the basis functions are real-valued, the coefficients in the expansion are quaternions (see[4]) .

For $f\ \in $ $L^1({\mathbb R})$ we denote by $\hat{f}$ , $\overline{\hat{f}\ }$the classical Fourier transforms of $f$~:  defined by   $\hat{f} \left(y\right)=\int^{+\infty }_{-\infty }{e^{-2\pi ixy}}\ f\left(x\right)dx$, \ $\overline{\widehat{f}}(y)=\frac{1}{\sqrt{2\pi }}\int^{+\infty }_{-\infty }{e^{-ixy}}\ f\left(x\right)dx.$

\begin{Lemma}

\[\widehat{{\varphi }_n}\ ={(-1)}^n{\varphi }_n.\] 
\end{Lemma}
 Proof. By using the functions $h_n\left(x\right)=\frac{{(-1)}^n}{n!}{{  e}}^{\frac{x^2}{2}}\frac{d^n}{{dx}^n}$(${{  e}}^{{  -}x^2}$)\ (see [3,  Proposition 1.9])
We get that  $\overline{\widehat{h_n}\ }\left(x\right)={(-1)}^nh_n\left(x\right),$

However ${\varphi }_n\left(x\right)=h_n\left(\sqrt{2\pi }x\right),$

Therefore $\widehat{{\varphi }_n}\ \left(x\right)$=$\int^{+\infty }_{-\infty }{e^{-2\pi ixy}}\ h_n\left(\sqrt{2\pi }y\right)dy$

\hspace*{3 cm}=$\frac{1}{\sqrt{2\pi }}\int^{+\infty }_{-\infty }{e^{-\sqrt{2\pi }ixz}}\ h_n\left(z\right)dz$

\hspace*{3 cm}=$\overline{\widehat{h_n}\ }$($\sqrt{2\pi }x$)

\hspace*{3 cm}=${(-1)}^nh_n\left(\sqrt{2\pi }x\right)$

\hspace*{3 cm}=${(-1)}^n{\varphi }_n\left(x\right).$

\begin{Lemma}

\[{\mathcal F}\left\{{\varphi }_{k,l}\right\}\left(y\right)={(-i)}^k{(-j)}^l{\varphi }_{k,l}\left(y\right).\] 
\end{Lemma}

Proof. \[{\mathcal F}\left\{{\varphi }_{k,l}\left(x\right)\right\}\left(y\right)=\int_{{{\mathbb R}}^2}{e^{-2\pi {  \ i}x_1y_1\ }{\varphi }_k\left({\ x}_1\right){\varphi }_l\left({\ x}_2\right)e^{-2\pi {  \ j}{\ x}_2y_2}dx_1dx_2}\] 

\hspace*{4cm}=$\int_{{\mathbb R}}{e^{-2\pi {  \ i}x_1y_1\ }{\varphi }_k\left({\ x}_1\right)} dx_1  \ \int_{{\mathbb R}}{{\varphi }_l\left({\ x}_2\right)e^{-2\pi {  \ j}{\ x}_2y_2}dx_2}$( Fubini)

\hspace*{4 cm}=${(-i)}^k{\varphi }_k\left({\ y}_1\right){(-j)}^k{\varphi }_l\left({\ y}_2\right)$ (lemma 5.1)

\hspace*{4 cm}=${(-i)}^k{(-j)}^l{\varphi }_{k,l}\left(y\right).$

\begin{Th}

Let $\alpha \ $  and $\beta $  are positive constants .Suppose $f\in L({\mathbb R}^2,{\mathbb H})$  with 

\[{|f\left(x\right)|}_Q\le {Ce}^{-\alpha {\left|x\right|}^2},  x\in {\mathbb R}^2.\] \hfill(5.3)

\[{|\ {\mathcal F}\left\{f\right\}\left(y\right)|}_Q\le {C\textprime e}^{-\beta {\left|y\right|}^2},\ y\in {\mathbb R}^2.\] \hfill(5.4)

for some positive constants $C,C\textprime.$Then, three cases can occur :

\textit{i)}  If $\alpha \beta >{\pi }^2$, then $f{  =0\ }$. 

ii ) If $\alpha \beta ={\pi }^2$, then $\ f{  =A}{{  e}}^{{  -}\alpha {\left|x\right|}^2}$ ,where A is a constant.

iii) If  $\alpha \beta <{\pi }^2,$\ then there many infinitely such functions $f$.
\end{Th}

Proof. We first prove the result for the case $\alpha \beta ={\pi }^2$.

With scaling we may assume $\alpha {  =}\beta =\pi $.

Indeed:

Let $g\left(x\right){  =}f\left(\sqrt{\frac{\pi }{\alpha }}x\right)$,then ${\mathcal F}\left\{g\right\}\left(y\right)=\frac{\alpha }{\pi }\ {\mathcal F}\left\{f\right\}\left(\sqrt{\frac{\alpha }{\pi }}y\right),$

Thus

\[{|{\mathcal F}\left\{g\right\}\left(y\right)|}_Q\le \frac{\alpha }{\pi }C\textprime e^{-\beta \frac{\alpha }{\pi }{\left|y\right|}^2}=\frac{\alpha }{\pi }C\textprime e^{-\pi {\left|y\right|}^2},\] 

and ${|g\left(x\right)|}_Q\le {Ce}^{-\pi {\left|x\right|}^2}.$

If the result is shown for $\alpha {  =}\beta =\pi $  

Then $g\left(x\right){  =A}{{  e}}^{{  -}\pi {\left|x\right|}^2}$ thus $f\left(x\right){  =}g\left(\sqrt{\frac{\alpha }{\pi }}x\right){  =A}{{  e}}^{{  -}\alpha {\left|x\right|}^2}.$

Now we assume \ $\alpha {  =}\beta =\pi $:

Complexifying the variable$\ \ y=a+i_{\mathbb C}{  \ b\ }$;${  \ a}$ =($a_1,a_2$),${  \ b}$ =($b_1,b_2$)$\ \in {\mathbb R}^2$ (we note by $i_{\mathbb C}$ the complex number checking $i^2_{\mathbb C}$= -1)

We have 

\[{\mathcal F}\left\{f\right\}\left(y\right)=\int_{{{\mathbb R}}^2}{e^{-2\pi {  \ i}x_1{  \ (}a_1{  +}i_{\mathbb C}\ b_1{  )}}f\left(x\right)e^{-2\pi {  \ j}{\ x}_2{  (}a_2{  +}i_{\mathbb C}\ b_2{  )}}dx_1dx_2},\] 

then

\[{|{\mathcal F}\left\{f\right\}\left(y\right)|}_Q\le \int_{{{\mathbb R}}^2}{{|f\left(x\right)|}_Q}e^{2\pi {  (}\left|x_1a_1\right|+\left|x_1b_1\right|{  +}\left|x_2a_2\right|+\left|x_2b_2\right|}dx_1dx_2\] 

According to the hypothesis ${|f\left(x\right)|}_Q\le {Ce}^{-\alpha {\left|x\right|}^2}.$

So

${|{\mathcal F}\left\{f\right\}\left(y\right)|}_Q\le Ce^{\pi {  (}{\left|a_1\right|^2+\left|b_1\right|^2)}}e^{\pi {  (}{\left|a_2\right|^2+\left|b_2\right|^2)}}\int_{{{\mathbb R}}^2}{e^{-\pi {  (}{\left|x_1\right|-(\left|a_1\right|+\left|b_1\right|))}^2}}e^{-\pi {  (}{\left|x_2\right|-(\left|a_2\right|+\left|b_2\right|))}^2}dx_1dx_2$

\hspace*{2.7cm}= $Ce^{\pi {\left|y\right|}^2}\int_{{\mathbb R}}{e^{-\pi {  (}{\left|x_1\right|-(\left|a_1\right|+\left|b_1\right|))}^2}}dx_1\int_{{\mathbb R}}{e^{-\pi {  (}{\left|x_2\right|-(\left|a_2\right|+\left|b_2\right|))}^2}}dx_2.$\\
Since $\int^{+\infty }_{-\infty }{e^{-\pi {  (}{\left|t\right|+m)}^2}}\ dt=\int^{+\infty }_0{e^{-\pi {  (}{t+m)}^2}}\ dt$ +$\int^0_{-\infty }{e^{-\pi {  (}{-t+m)}^2}}\ dt\ $ for $m\in {\mathbb R}$

\hspace*{4cm}=$\int^{+\infty }_0{e^{-\pi {  (}{t+m)}^2}}\ dt$+$\int^{+\infty }_0{e^{-\pi ({t+m)}^2}}\ dt$\\
\hspace*{4.4cm}=2$\int^{+\infty }_m{e^{-\pi t^2}} dt$\\
\hspace*{4.4cm}  $ \le 2 \int^{+\infty }_{-\infty }{e^{-\pi t^2}}\ dt=\ 2.$

We deduce that ${|{\mathcal F}\left\{f\right\}\left(y\right)|}_Q\le 4 \  Ce^{\pi {\left|y\right|}^2}.$\hfill(5.5)

As ${\mathcal F}\left\{f\right\}\ $is an entire function satisfying relation (5.4)and (5.5), lemma 2.1 in [10]

Allows to express ${\mathcal F}\left\{f\right\}\ $as follows :

${\mathcal F}\left\{f\right\}\left(y\right)=Ae^{-\pi {\left|y\right|}^2}$ with A is constant.

Since ${\mathcal F}\left\{e^{{  -}\pi {\left|x\right|}^2}\right\}\left(y\right)$=$e^{{  -}\pi {\left|y\right|}^2}$(by Lemma 3.5), the injectivity of the two sided quaternion Fourier transform(Lemma 3.2) completes the proof.

For the case $\alpha \beta >{\pi }^2$.

 By  scaling we may assume  $\alpha =\beta >\pi \ $then

${|{\mathcal F}\left\{g\right\}\left(y\right)|}_Q\le C\' e^{-\pi {\left|y\right|}^2}$and ${|g\left(x\right)|}_Q\le {Ce}^{-\pi {\left|x\right|}^2}$ for some constants $C,C'>0.$

Then according to the first case ,$\ f\left(x\right){  = A}{{  e}}^{{  -}\pi {\left|x\right|}^2}$,where A is a constant, 

but this contradicts the bound \textit{ }${|f\left(x\right)|}_Q\le {Ce}^{-\alpha {\left|x\right|}^2}$unless\textit{ }$f=0.$

For the final case $\alpha \beta <{\pi }^2$

We again assume $\alpha = \beta <\pi, $

we show by induction (see  [9, Thm. 2, p. 5] ) \ that ${|{\varphi }_k\left(x\right)|}_Q\le {Ce}^{-\alpha x^2}$ for $x\in {\mathbb R},$ $\alpha <\pi $   (${{\{\varphi }_k\}}_{k \in {\mathbb N}}$ are real-valued)

Thus ${|{\varphi }_{k,l}\left(x\right)|}_Q$=${|{\varphi }_k\left(x_1\right)|}_Q{|{\varphi }_l\left(x_2\right)|}_Q \ \le {C'e}^{-\alpha {\left|x\right|}^2}.$ 

Now we need only show that ${\mathcal F}\left\{{\varphi }_{k,l}\right\}$\ obeys similar bound.

We have

${|{\mathcal F}\left\{{\varphi }_{k,l}\left(x\right)\right\}\left(y\right)|}_Q$=${|{(-i)}^k{(-j)}^l{\varphi }_{k,l}\left(y\right)|}_Q$(Lemma 5.2)

\hspace*{3.3 cm}=${|{\varphi }_{k,l}\left(y\right)|}_Q$ \\
\hspace*{4cm}
$({|{(-i)}^k{(-j)}^l{\varphi }_{k,l}\left(y\right)|}_Q={|{(-i)}^k|}_Q{|{(-j)}^l|}_Q{|{\varphi }_{k,l}\left(y\right)|}_Q={|{\varphi }_{k,l}\left(y\right)|}_Q$)

\hspace*{3.3 cm}$\le {Ce}^{-\alpha {\left|y\right|}^2}$.\\
This completes the proof of the Hardy's theorem.


\end{document}